\input amstex
 \documentstyle{amsppt}
 \magnification=1200
 %\NoRunningHeads
 \vsize=8.5truein
 \hsize= 6.0truein
 \hoffset=0.3truein
 \voffset=0.2truein
  \TagsOnRight
\define\db{\bar\partial}

\nologo
 \NoRunningHeads
 \NoBlackBoxes

 \topmatter
\title Pseudo-Einstein and Q-flat metrics with eigenvalue estimates on
CR-hypersurfaces
\endtitle
\author Jianguo Cao\footnote{Supported partially by NSF Grant
DMS 0404558. The first named author is very grateful to National
Center for Theoretical Sciences at National Tsinghua University
for its warm hospitality. \hfill{$\,$}} and Shu-Cheng Chang
\footnote{Supported partially by the NSC of Taiwan \hfill{$\,$}}
\endauthor
\address Mathematics Department, University of Notre Dame,
Notre Dame, IN 46556
\endaddress
\email jcao\@nd.edu
\endemail
\address Department of Mathematics, National Tsing Hua University,
Hsinchu 30013, Taiwan, R. O. C.
\endaddress
\email scchang\@math.nthu.edu.tw
\endemail

\abstract

    In this paper, we will use  the Kohn's $\db_b$-theory on
    CR-hypersurfaces to derive some new results in CR-geometry.

\proclaim{Main Theorem} Let $M^{2n-1}$ be the smooth boundary of a bounded strongly pseudo-convex domain $\Omega$ in a complete Stein manifold
$V^{2n}$.   Then (1) For $n \ge 3$, $M^{2n-1}$ admits a pseudo-Einstein metric; (2) For $n \ge 2$,  $M^{2n-1}$ admits a Fefferman  metric of zero
CR Q-curvature; and (3) for a compact strictly pseudoconvex CR emendable 3-manifold $M^3$, its CR Paneitz operator
$P$ is a closed operator.

\endproclaim

There are examples of non-emendable  strongly pseudoconvex CR
manifold $M^3$, for which the corresponding $\db_b$-operator and
Paneitz operators are not closed operators.
\endabstract
\endtopmatter

\document
 \baselineskip=.25in

     \vskip2mm

\vskip2mm \vskip2mm

 \head 0. Introduction
 \endhead

    \vskip4mm

    In this paper, we study several questions, including the existence of  $Q$-flat metrics, pseudo-Einstein metrics and
    the closedness of the CR Paneitz operators.

    First, we will use an approach proposed by Fefferman and his school to prove that  ``the complete K\"ahler-Einstein $g_\infty$
    on an open domain $\Omega$ induces a metric on $M = b\Omega$ with zero CR $Q$-curvature, where
    $\Omega$ is a smooth, bounded
    strictly pseudo-convex domain in a Stein manifold $V^{2n}$."    To achieve this goal, we solve a $\partial \db$-Poincar\'e-LeLong
    equation
    via the $\db$-theory. Although this part does not produce new hard a-priori estimates, it is still valuable for other potential applications.

    The second purpose is to prove the existence of pseudo-Einstein metrics on strictly pseudo-convex
    $CR$-hypersurface of real dimension $\ge 5$ through solving the $\db_b$ Poincar\'e-LeLong equations.

       The last part of our paper is to study the closedness of CR Paneitz operator, which is a fourth-order differential
    operator. It is known that the positivity of CR Paneitz operator is related to the deformation of  $Q$-curvatures  under
    the conformal change of metrics on Riemannian manifold $M^m$. In particular, the positivity of
    CR Paneitz operator is also related to
    the lower bound of the first eigenvalue of sub-Laplace on a $CR$ manifold $M^{3}$, see [CC], [CCC] and [LL].
    It will be shown that, if $M^3 = b\Omega^4$ is the smooth
    boundary of bounded strictly pseudo-convex domain  $ \Omega$ in a Stein manifold $V^4$, then its CR-Paneitz
    operator on $M^3$ is closed, for any metric on $M^3$.

\proclaim{Main Theorem} Let $M^{2n -1}$ be the smooth boundary of a bounded strongly pseudo-convex domain $\Omega$ in a complete Stein manifold
$V^{2n}$.   Then

\smallskip
\noindent
 (1) For $n \ge 2$,  $M^{2n -1}$ admits a metric of zero CR Q-curvature;

 \smallskip
\noindent
 (2)  For $n \ge 3$, $M^{2n - 1}$ admits a pseudo-Einstein metric;

 \smallskip
\noindent
 (3) In addition, for a compact strictly pseudoconvex CR emendable
 3-manifold $M^3$, its CR Paneitz operator
$P$ is a closed operator.

\endproclaim

Earlier work in this direction for the case of $V^{2n} = \Bbb C^{n}$ can be found in [L2], [FH] and [GG]. In a very recent paper [LL], Li and
Luk obtained an explicit formula for Webster's pseudo-Ricci curvature on real hypersurfaces in $\Bbb C^n$.  Thus, their result could lead another
proof of Cheng-Yau's result ([CY]) and Mok-Yau's theorem [MY], which will be used in Section 2 below.

Among other things, we introduce some new methods to handle pseudo-Einstein metric and Paneitz operators in this paper. For example, we use the
closeness of $\db_b$ and $\db_b^*$ operators provided
 by Kohn's theory, in order to
complete the proof. When $\dim_{\Bbb R} [M] = 3$, we decompose the Paneitz operator $P$ as a product of closed operators. Thus, the
closed property of $P$ will follow immediately, see Lemma 1.4 and Section 4 below.

\vskip2mm \vskip2mm

 \head 1. Preliminary results
 \endhead

    \vskip4mm
    It is well-known that the real Laplace $\triangle$ on a
    K\"ahler manifold $V^{2n}$ satisfies
    $$
\triangle = 2 \square = 2 \bar \square,
    $$
where $\square = \db \db^* + \db^* \db$ is the complex Laplace operator. However, it may happen that $\triangle_b  \neq  2 \square_b$ in some cases.
Let us recall the notions  of $\triangle_b $ and $ \square_b$.

Since $M^{2n-1} = b V^{2n}$ has odd real dimension, it is a Cauchy-Riemann manifold. The $\db_b$
operator induces
a sub-elliptic operator
$$
\square_b = \db_b^* \db_b + \db_b \db_b^*
$$
acting on $L^2_{(p, q)}(M)$. Similarly, there is a real sub-Laplace operator, which can be viewed as partial trace of the hessian operator (or
can be viewed a sum of the squares of $(2n-2)$ vectors):
$$
       \triangle_b u |_z = \sum_{k = 1}^{2n-2} \langle \nabla_{e_k} (\nabla^b u), e_k \rangle |_z
$$
where $e_{2n}$ is the outward real unit normal vector of $\Omega$ along boundary $M = b \Omega$, $e_{2j} = J e_{2j-1}$ for
$j = 1, ..., n$, $z \in M$, $J $ is the complex structure of $V^{2n}$,
$\{e_1, e_2, \cdots, e_{2n-3}, e_{2n-2},
e_{2n-1}, e_{2n} \}$ is an orthonormal basis of $[T_z(V)]_{\Bbb R}$ and
$$
\nabla^b u = \sum_{k = 1}^{2n-2} du( e_k) e_k.
$$

When $\Omega$ has the strongly pseudo-convex boundary in a Stein manifold $V^{2n}$ with $n =2$, it has been observed that
$$
   \square_b u = \frac 12  [\triangle_b u + \sqrt{-1} T u] \tag1.1
$$
for all $u \in L^2(M^3)$, where $T = \lambda e_3 $ is the Reeb vector of the CR 3-manifold $M^3$ for some real valued function $\lambda$, see
[L1, p414].

    The operator $\square_b$ is a Lewy type operator, which may {\it not} be locally solvable.

    If the Reeb vector $T$ induces an infinitesimal pseudo-conformal with respect to the Tanaka-Webster metric,
    then the torsion of $M^3$ is zero, see [Web, p33]. In this case, the operator $\square_b$ is related
    to the so-called CR Paneitz operator $P$, where $P$ is given by
    $$
P u = \triangle_b^2 u +  T^2 u = 4 \square_b \bar{\square}_b   u, \tag1.2
    $$
for $u \in L^2(M^3)$.
More generally, if $M^3$ has torsion free in the sense of Tanaka (cf. [Ta1-2] [Web]), then (1.2) holds.

    The eigenvalues of the Paneitz operator and CR Paneitz operators have been considered various authors
    ([Ch], [CC]). The eigenvalue
    estimate plays an important role
    to the study of the so-called Q-curvature flow, see [Br] [CCC].

\proclaim{Definition 1.1} (1)  The CR-Paneitz operator $P: L^2(M^3) \to L^2(M^3)$ is called essentially positive, if there is a
positive constant $\lambda_1 >0 $ such that
$$
  \langle  P u, u \rangle \| \ge \lambda_1 \| u \|^2,  \tag1.3
$$
for all $ u \bot ker(P)$.

(2) The operator $\Cal F: L^2_{(p, q)}(M) \to L^2_{(p, q)}(M)$ is said to have
positive spectrum gap at $0$  (or is said to be a closed operator) if
there is a
positive constant $\lambda_{p, q} >0 $ such that
$$
  \| \Cal F u \| \ge \lambda_{p, q} \| u \|,  \tag1.4
$$
for all $ u \bot [L^2_{(p, q)}(M) \cap ker(\Cal F)]$.

(3) A smooth function $ f: U_\varepsilon(M) \to \Bbb R$ is called a defining function of $M$ if $f^{-1}(0) = M$ and if
$0 $ is not a critical value of $f$, where $U_\varepsilon(M) \subsetneq V^{2n}$ is a neighborhood of $M$ in a Stein manifold $V^{2n}$.

(4) Let $\theta$ be a contact 1-form of $M^{2n-1}$ and $J: \ker
\theta \to \ker \theta$ be the almost complex structure on the
CR-distribution $\ker \theta$ such that $J^2 \vec v = - \vec v$
for all $\vec v \in \ker \theta$. In what follows, we always let
$$
[T^{(1, 0)}(M) \oplus T^{(0, 1)}(M)] = [\ker \theta] \bigotimes_{\Bbb R} \Bbb C.
$$

(5) A CR manifold $M^{2n-1}$ is said to have transverse symmetry or
torsion-free if it admits a CR Reeb vector field $\xi $ such that $\xi \notin \ker \theta$ with
$$
\Cal L_\xi J = 0
$$
where $\Cal L$ is the Lie derivative and $J$ is the complex structure of $[T^{(1, 0)}(M) \oplus T^{(0, 1)}(M)]$.
\endproclaim

If $\xi $ is the real part of a holomorphic vector filed $\tilde X$ on a neighborhood  $U_\varepsilon(M)$ of $M$, then $\xi$ induces an automorphism on
$U_\varepsilon(M)$. Any real
part $\xi$ of a holomorphic vector filed restricted to $M$ induces a CR-automorphism of $M$.

   In the H\"ormand-Kohn $L^2$-theory and the Kohn-Rossi theory, the essential spectrum of $\square$ and $\square_b$ have
   been extensively investigated.

   A smooth $(p, q)$-form  $u$ on $\Omega$ with $q \ge 1$ is said to satisfy the $\db$-Neumann boundary condition if
   $$
           u( (\db \rho)_\#, ...) |_z = 0
   $$
   for all $z \in M = b\Omega$, where $(\db \rho)_\# $ is the complex normal vector field of type $(0, 1)$
   along the boundary  $M^{2n-1}$.

\proclaim{Theorem 1.2} ([CS], [CaWS]) Let $\Omega$ be a bounded domain
with smooth pseudo-convex boundary $M$ in a complete Hermitian
manifold $V^{2n}$. Suppose that $V^{2n}$ is either a Stein
manifold or $\Bbb CP^n$. Then  the complex Laplace
operator $\square$ is

\smallskip
\noindent
(1) positive for on $L^2_{(p, q)}(\Omega)$ with $(n-1) \ge q\ge 1$; and

 \smallskip
\noindent
(2)  essentially positive on $L^2_{(p, 0)}(\Omega)$ and $L^2_{(p, n)}(\Omega)$

\smallskip
\noindent
with respect to
$\db$-Neumann boundary condition on $M = b\Omega$.

Moreover, for any Hermitian metric on $\Omega$, the operator $\square$ is essentially positive
on $\Omega$ with respect to $\db$-Neumann
boundary condition on $M$.
\endproclaim

For the $L^2$ estimates of $\square$, the domains $\Omega$ in
Theorem 1.2 are not necessarily strictly pseudo-convex.
However, for estimates of $\square_b$ on the boundary $M^{2n-1}$ of $\Omega$, we need extra assumptions on $M^{2n-1}$.

The dual of $\db $-Neumann problem is the so-called $\db$-Cauchy
problem. A $(p, q)$-form $u$ is said to satisfy the Cauchy
boundary condition on $M = b\Omega$ if
$$
        u( \xi, ...)|_z = 0
$$
for all $\xi \in T_z^{(0, 1)}(M)$ and $z \in M$. If a $\db$-closed form $f \in C^\infty_{(p, q+1)}(\Omega)$ with a compact support in $\Omega$, then one
consider to solve $\db u = f$ such that $u$ has a compact support in $ \Omega$ as well. Solving $\db u = f$ with compact support is related to the
$\db$-extension problem, via the Kohn-Rossi theory. Using the solution to the $\db$-extension problem and Theorem 1.2, we are able to solve
$\db_b u = f$ on a special class of CR-manifolds:

 \proclaim{Theorem 1.3}  ([CS], [CaSW])  Let $\Omega$ be a bounded Hermitian manifold with a smooth pseudo-convex boundary $M$.
Suppose that one of the following
 conditions holds:

(1) $\Omega$ is a domain of a  complete Stein manifold $  V^{2n}$;

(2) $\Omega \subset \Bbb CP^n$, and  $M = b\Omega$ admits a pluri-subharmonic defining function.

Then
 the $\db$-Cauchy boundary problem is solvable on $\Omega$. Furthermore, (1) $\db_b$-operator is closed; and
   (2) the operator $\square_b: L^2_{(p, q)}(M) \to L^2_{(p, q)}(M)$ is  positive for $1 \le q \le n-2$ and essentially positive
   for $q =0$ or $q = n-1$.
\endproclaim

When $M = b\Omega$ is strongly pseudo-convex, it is well-known that $M$ admits a pluri-subharmonic defining function, see [DF].

If $\Cal L: H_1 \to H_2$ is a linear operator, we let $\text{Dom}(\Cal L)$ be its domain and $\Cal R (\Cal L)   $ be its range. If $ A \subset
H$ is a subset of a Hilbert space $H$, the closure of $A$ in $H$ is denoted by $\bar{A}$.

We begin with an elementary but useful criterion for closed operators.

\proclaim{ Lemma 1.4} ([CS, p60] or [H\"o1-2]) Let $\Cal L: H_1 \to H_2$ be a linear, closed, densely defined operator from the Hilbert space $H_1$
to another Hilbert space $H_2$. The following conditions on $\Cal L$ are equivalent:

\smallskip
\noindent (1) The range $\Cal R (\Cal L)   $ of $\Cal L$ is closed;

\smallskip
\noindent (2) There is a constant $C$ such that
$$
             \| f \|_1  \le C \| \Cal L f \|_2
$$
for all $f \in \text{Dom}(\Cal L) \cap { \Cal R (\Cal L^*)}$;

\smallskip
\noindent (3)  The range $\Cal R (\Cal L^*)   $ of $\Cal L^*$ is closed;

\smallskip
\noindent (4) There is a constant $C$ such that
$$
             \| f \|_2  \le C \| \Cal L^* f \|_1
$$
for all $f \in \text{Dom}(\Cal L^*) \cap \Cal R (\Cal L)$.
\endproclaim

\head 2. The existence of CR $Q$-flat metrics  on strictly pseudo-convex CR-hypersurfaces  in a Stein manifold
\endhead
\vskip2mm

In this section, we first recall an existence result of CR $Q$-flat metrics on CR-hypersurfaces in Euclidean space $\Bbb C^n$ due to Fefferman
and others. Afterwards, we will extend such a
result to $CR$-hypersurfaces in an arbitrary Stein Manifold $V^{2n}$.
One of our key steps is to use the $\db$-theory to introduce the generalized Fefferman's functional
$u \to \hat{J}(u)$, which is independent of the choice of local holomorphic coordinates, see (2.5) below.

\subhead 2.a. A sufficient condition for existence of  $Q$-flat metrics
on real hypersurfaces
\endsubhead

    Let us recall a sufficient condition for existence
   of $Q$-flat metrics on real hypersurfaces, which were derived by Fefferman and others.

   \proclaim{Proposition 2.0  } ([FG1-2], [GG]) Let $\Omega \subset \Bbb C^n $ be a compact domain
   with smooth boundary $M^{2n-1} = b \Omega$ in the complex Euclidean space $\Bbb C^n$.
   Suppose $\Sigma^{2n}$ is an unit circle bundle defined on a
   $CR$-hypersurface $M^{2n-1}$ and suppose that $\Sigma^{2n}$  admits an $S^1$-invariant
   Einstein-Lorentz metric $g_u^+ = i \partial \db H_u |_{\Sigma^{2n}} $ defined as below. Then $M^{2n-1}$ admits a metric of zero
  $CR$ $Q$-metric.
\endproclaim

    We now provide a description of the metric $g_u^+$ stated in Proposition 2.0, which will be used for
    any real hypersurface $M^{2n-1} $ in a Stein manifold $V^{2n}$ as well.

   Let $K^*$ be the canonical bundle of $V^{2n}$ restrict to $M$ and let $\Sigma^{2n} = K^*/ {\Bbb R^+}$ be
   the unit circle bundle of $K^*$. Thus there is
a fiberation
$$
         S^1 \to \Sigma^{2n} \to M^{2n -1}
$$
and $\dim_{\Bbb R}(\Sigma^{2n}) = 2n $.

   We may assume that $\Omega \subset V^{2n}$ is an open strictly pseudo-convex domain with  compact smooth boundary $M^{2n - 1} = b\Omega$.
   Suppose that $\hat u $ is a defining function of $M^{2n - 1}$. For example, we can choose $\hat u $ as a signed distance function form $M$:
   $$
\hat u(z) = \biggl\{ \aligned  & -d(z, M),  \quad if \quad z \in \Omega \\
& d(z, M), \quad if \quad z \notin \Omega
\endaligned
   $$
Any other defining function $u $ can be expressed as
$$
   u(z) = e^\eta \hat u
$$
for some real valued function $\eta$.

The contact structure on $M$ is an 1-form given by
$$
  \theta_u (\xi) = du ( J \xi)
$$
for all $\xi \in [T(M)]_{\Bbb R}$, where $J$ is the complex structure of $V^{2n}$.

    There are two types of metrics which we will use. The first one is the Cheng-Yau metric on $\Omega$; and the second one is
    introduced by Fefferman on a line bundle over $b\Omega$.

     Let us first consider complete K\"ahler metrics on an open domain $\Omega$. Suppose
that
$$
          \omega_u = i \partial \db [ \log ( - \frac 1u     )]
$$
is a K\"ahler form on $\Omega$. Such a K\"ahler form $\omega_u$ corresponds to a K\"ahler metric
$$
   g_u ( X, Y) = \omega_u (X, JY) = i \partial \db [ \log (  -\frac 1u     )](X, JY), \tag2.1
$$
where $J$ is the complex structure of $\Omega$.

Secondly, Fefferman and his school considered a class of Lorentz metrics on
     canonical bundle on $K^*$ mentioned above.

     We will use an extrinsic way to define such metrics, along the line described in a new book [DT, p150].
Suppose that $\Lambda_{(n, 0)}(V^{2n}) $ be the canonical line bundle of open domain $V^{2n}$. Clearly,
$  \Cal L_{V^{2n}}= \Lambda_{(n, 0)}(V^{2n}) $ is a complex manifold of complex dimension $(n+1)$.

When $\xi$ is a cross-section of $\Cal L_{V^{2n}}$ over $V^{2n}$,
the norm $|\xi|_{g_u}$ induced by $g_u$ is well-defined. We
further define
$$
       H_u(z, \xi) = |\xi|_{g_u}^{\frac{2}{n+1}}u(z)
$$
There is an $(1,1)$-form defined on $\Cal L_{V^{2n}}$ given by
$
             i \partial \db H_u.
$

Similarly, there is a Hermitian form
$$
           G_u ( \tilde X, \tilde Y) =  i \partial \db H_u (\tilde X, \tilde J \tilde Y), \tag2.2
$$
where $\tilde J$ is the complex structure of line bundle $\Cal L_{V^{2n}}$. The Hermitian form $G_u$ is
{\it not}  necessarily positive definite on
the complex manifold $\Cal L_{V^{2n}}$.

We now consider a subset
$$
 \Sigma^{2n} = \{ (z, \xi) \in \Cal L_{V^{2n}} \quad | \quad z \in b\Omega, |\xi| = 1  \} \tag2.3
$$
where $\Omega$ is an open, bounded and strictly pseudo-convex domain in $V^{2n}$.

Finally, when $i \partial \db u > 0 $ on $M = b\Omega$,  we consider
$$
g^+_u = G_u |_{\Sigma^{2n}}.    \tag2.4
$$
It was shown that $g^+_u $ is a Lorentz metric on $\Sigma^{2n}$. Clearly, $\Sigma^{2n}$ is diffeomorphic to the unit circle bundle $K^*$ mentioned above.

We remark that  the function $u =0$ {\it vanishes} on $M^{2n-1}$. The leading term of the metric $g^+_u$ is
$$
 i \partial \db u.
$$

In [FH], Fefferman and Hirachi studied the so-called $Q$-curvature of $CR$-manifold $M^3$:
$$
    Q^{CR}_{\theta_u} = \frac 43 (\Delta_b R - 2 Im \nabla^\alpha \nabla^\beta A_{ \alpha\beta}),
$$
where $R$ is the Tanaka-Webster scalar curvature, $A$ is the torsion, $\Delta_b$ is the sub-Laplacian computed in terms of the
contact 1-form $\theta_u$ and $\theta_u(\xi) = d u( J \xi)$ for all $\xi \in T(M)$.

For higher dimensional manifolds,  the $Q$-curvatures of higher order have been studied in [FH] and [GG].

The notations above will be used in the next two sub-sections.

\medskip

\subhead 2.b. Relations between the Fefferman's Lorentz metric and the
Cheng-Yau's K\"ahler-Einstein metric
\endsubhead

In this sub-section, we illustrate a strategy to obtain the existence of $Q$-flat metrics on real hypersurfaces in $\Bbb C^n$.

Let us now recall a result obtained by Fefferman and his school.

\proclaim{Proposition 2.1}  ([FG1, Chapter III]) Let $\Omega \subset V^{2n}$, $M = b\Omega \subset \Bbb C^n$, $u = \hat u e^\eta$ and $\{g_u,  g^+_u \}$ be as above.
 If the Cheng-Yau metric
$g_u$ is a complete K\"{a}hler-Einstein on $\Omega$, then the Lorentz metric
$g^+_u $ is Einstein on $\Sigma^{2n}$.

\endproclaim

Here is a direct application of Propositions 2.0-2.1.

\proclaim{Corollary 2.2}  ([FH], [GG]) Let $\Omega \subset \Bbb C^{2n}$ be an open strictly pseudo-convex
domain with compact closure  and let $M^{2n-1} = b\Omega$ be its boundary. Then
$M$ admits a metric of zero CR $Q$-curvature.
\endproclaim

Proposition 2.1 and Corollary 2.2 were stated for strictly pseudo-convex and bounded
   domain $\Omega$ in $\Bbb C^n$. We would like to extend these results to
   any strictly pseudo-convex and bounded
   domain $\Omega$ in  a Stein manifold $V^{2n}$.

\medskip

\subhead
2.c. Compact smooth real hypersurfaces in a Stein manifold
\endsubhead

\medskip

   Our goal of this section is to verify the following theorem.

         \smallskip
\proclaim { Proposition  2.3} Let $\Omega  $ be a bounded, open and strictly pseudo-convex domain with
a smooth boundary in a Stein manifold $V^{2n}$. If the metric $g_u$
above is a complete K\"ahler-Einstein metric on $\Omega$, then $ g_u$ induces a
metric $\tilde{ g}^\infty_u$ on $M = b\Omega$ with zero CR $Q$-curvature.
\endproclaim
\demo{Proof} Since $V^{2n}$ is Stein, we may assume that $V^{2n} \subset \Bbb C^m$ is a complete sub-manifold of $\Bbb C^m$, for sufficiently large $m$.
   Let $\hat g $ be induced metric on $\Omega \subset V^{2n} \subset  \Bbb C^m$. For each local holomorphic coordinate
   system $\{(z_1, ..., z_n) \}$ of $\Omega$, the Ricci tensor $\hat Ric $ of $\hat g$ is given by
   $$
           \hat Ric = -  i \partial \db \log [\det \hat g_{i \bar j}].
   $$
It is clear that $\hat Ric$ is well-defined and independent of the choice of local holomorphic coordinate
   system $\{(z_1, ..., z_n) \}$. Moreover, $\hat Ric$ is a closed $(1,1)$-form on $\Omega$.
In what follows, we first would like to solve Poincare-Lelong equation $i \partial \db  f =   \hat Ric$.

For this purpose, we recall a theorem of Dolbeault:
 $$ H^{(1, 1)}(\Omega) = H^{(0, 1)}(\Omega, \Cal O|_\Omega)
 $$
where $\Cal O|_\Omega$ is the bundle of holomorphic $(1, 0)$-forms.

      Since $\Omega$ is strictly pseudo-convex and bounded
   domain in a Stein manifold $V^{2n}$, by  a theorem of Andreotti and Vesentini [AV], we have
   $$
 H^{(0, 1)}(\Omega, \Cal O|_\Omega) = 0.
   $$
In fact, Proposition A.4 of [CaWS, p218] is also applicable for (0, q)-forms with values  in $ \Cal O|_\Omega$.  Thus,
$H^{(1, 1)}(\Omega) = H^{(0, 1)}(\Omega, \Cal O|_\Omega) = 0$. Professor Siu  also handled similar formula with values in a
vector bundle $E$, although the weighted functions were not discussed there (cf. [Siu, Chapters 2-3]). Hence,
the first Chern class $c_1(\Cal O|_\Omega  ) = 0$. Recall that, by Chern-Weil theory, the co-homology class
$c_1(\Cal O|_\Omega  )$ is independent of the choices of affine connections, (cf. [Mi]). Therefore,
$c_1(\Cal O|_\Omega  ) = 0$
implies that
the Chern-Weil form $\hat Ric$
is $d$-exact on $\Omega$.

Therefore, we have $\hat Ric = d \beta$ for some $1$-form $\beta$.
Let us consider the decomposition of $\beta = \beta^{(1, 0)} + \beta^{(0,1)}$, where $\beta^{(0,1)}$ is the $(0, 1)$-component of $\beta$.
If $\hat Ric = d \beta$ and if $\beta = \beta^{(0, 1)} + \beta^{(1,0)}$, then $\db \beta^{(0, 1)} = 0$, where we used
the fact that $\hat Ric$ is an (1,1)-form. Choosing $f$ with $\db f = i \beta^{(0, 1)}$, we get a solution $i \partial \db  f =   \hat Ric$.

 Recall that $ \hat Ric$ is real valued. Replacing $f$ by $Re\{f\}$ if needed, we conclude that  the Poincare-Lelong equation
   $$
       i \partial \db  f =   \hat Ric = -  i \partial \db \log [\det \hat g_{i \bar j}].
   $$
has a smooth real-valued solution $f$ on $\Omega \cup b\Omega$. Such a solution $f$ is unique up to adding a pluri-subharmonic function. If we
require that $f$ has the smallest $L^2(\Omega)$-norm, then such a solution is unique, see Chapters 4-5 of [CS]. Such a solution $f$ is called a
Ricci potential of $\hat g$.

Following Fefferman [F2], we consider
$$
          \hat{J}(u) = (-1)^n e^{-f} \frac{1}{ \det \hat g_{i \bar j}} \det \left( \matrix  u & u_{\bar j} \\
                                            u_i & u_{i \bar {j} }  \endmatrix \right) \tag2.5
$$
where $f$ is the Ricci potential of $\hat g$ as above, $u_i = \frac{\partial u}{\partial z_j} $, $u_{i \bar {j} } = \frac{\partial^2 u}{\partial z_i \partial \bar{z}_j }$ and
     $\{z_1, ..., z_n \}$ is a local holomorphic frame.

     When $\Omega \subset \Bbb C^n$, we choose the standard coordinate system. Thus, in this case,
$ \det \hat g_{i \bar j} = 1$ and we can choose $f = 0$. Therefore, our definition coincides with
Fefferman's definition for the case of $\Omega \subset \Bbb C^n$, see [F2] and [CY].

     A calculation similar to [CY, p508] further shows that the metric $g_u$ is K\"ahler-Einstein of
     negative curvature
     $-(n+1)$  if
     $$
 \frac{\det \varphi_{i \bar j} }{ \det \hat g_{i \bar j}} = e^f e^{(n+1)\varphi} \tag2.6
     $$
holds, where $ \varphi =  \log (- \frac 1u)$.

A further calculation shows that the above equation holds if and only if
$$
      \hat{ J}(u)|_z \equiv 1  \tag2.7
$$
holds for all $z \in \Omega$.

 It is known that if $\hat{J}(u)|_z \equiv 1$ in $\Omega$, then $M = b\Omega$ has
zero CR Q-curvature, see [FH, Chapter 3]. This completes the proof of Proposition 2.3.
   \qed\enddemo

\proclaim{Corollary 2.4} Suppose that $\Omega \subset V^{2n}$ be a bounded, open and strictly pseudo-convex domain
with smooth boundary in
a Stein manifold $V^{2n}$ with $n \ge 2$. Then its boundary $M^{2n-1} = b\Omega$ admits a metric of zero $Q$-curvature.
\endproclaim
\demo{Proof} By Proposition 2.3,
   it remains to verify that there is a complete K\"ahler-Einstein metric $g_u$ on $\Omega$. The existence of
   such a complete K\"ahler-Einstein metric $g_u$ is provided by Mok-Yau in [MY, p52]. In fact, Mok and Yau found  desired solutions
   $u = e^\eta \hat u$ and $ \varphi =  \log (- \frac 1u)$ satisfying
   $ \frac{\det \varphi_{i \bar j} }{ \det \hat g_{i \bar j}} = e^f e^{(n+1)\varphi}    $. \qed\enddemo

\head 3. Existence of Pseudo-Einstein metrics on CR-hypersurfaces of real dimension $\ge 5$
\endhead

  In this section, we discuss the existence of pseudo-Einstein metrics on CR-hypersurfaces of real dimension $\ge 5$.
A metric $g$ defined on a $CR$-manifold $M^{2n-1}$ is said to be {\it pseudo-Einstein (or partially Einstein)}
if its Ricci tensor satisfies
$$
       Ric_g (X, Y)|_z = \lambda g(X, Y)|_z  \tag3.0
$$
for some constant $\lambda = \lambda(z) $ and for all real vectors $\{X, Y\}$ in the $CR$-distribution $ker(\theta)|_z$, where
$\theta $ is the contact form of $M^{2n-1}$.

  One of our new contributions in this section is to use the $\db_b$-theory to solve
  boundary version of Poincar\'e-Lelong equation related to the {\it partially Einstein equation}, see
  Proposition 3.4 and Corollary 3.5 below.

    When $\dim_{\Bbb R}[M^{2n-1}] = 3$, any metric $g$ on $M^3$ is {\it pseudo-Einstein (i.e., partially Einstein)}. Therefore, we
    only consider the case of $\dim_{\Bbb R}[M^{2n-1}] \ge 5$.

    We emphasize that a pseudo-Einstein metric $g$ on $M^{2n-1}$  is {\it not} necessarily Einstein. The pseudo-Einstein condition
    puts {\it no} restriction on its Ricci curvature in the directions which are transversal to $CR$-distribution. It might happen that
    $$
Ric_g (Z, Y) \neq  \lambda g(Z, Y)
    $$
    for some transversal vector $Z \bot ker( \theta )$.

  In [L2], Lee already showed that, if a compact strongly pseudo-convex CR-manifold $M^{2n-1}$ admits a closed, nowhere vanishing
  $(n, 0)$-form, then $M^{2n-1}$ admits a pseudo-Einstein metric. In particular, if $M = b \Omega$ and $\Omega \subset \Bbb C^n$, then
  $M$ admits a pseudo-Einstein structure.

  We make extra observations to extend Lee's  result to the case of $\Omega \subset V^{2n}$
  for any Stein manifold $V^{2n}$. The new ingredient of our approach will
  use the fact that the Chern curvature forms $\Theta$ are type of $(1, 1)$ for Lorentz-K\"ahler metrics.

  In addition, we will use Kohn's $\db_b$-theory to solve the boundary version of
 Poincare-Lelong equation
  $$
       i \partial_b  \db_b  f = \Theta \tag3.1
  $$
for any $\db_b$-closed $(1, 1)$-form $\Theta$.

The equation (3.1)
above is related to the existence of pseudo-Einstein metrics, as
described in [L2, p173].  Such an equation
was previously studied in [CaWS] for other purposes.

It is well-known that, for any function $u$, one has
$$
(d^c u) (\xi) = (du)  (J \xi) \text{  and  }  d d^c u = i \partial
\db u.
$$

   We begin with an elementary observation.

\proclaim{Lemma 3.1} Let $\hat u$ be a defining function of $M = b\Omega$. Suppose that
$\Omega \subset V^{2n}$ is a strictly pseudoconvex bounded domain in a Stein manifold.
Then

\smallskip
\noindent
(1) There is another defining function $u = e^\varphi \hat u $ such that $u$ is a strictly
pluri-subharmonic in a neighborhood of $M = b\Omega$, i.e., $i \partial \db u > 0$.

\smallskip
\noindent
(2) When $i \partial \db u > 0$ and $\theta_u = d^c u$, then $i \partial \db u$ gives rise to
a K\"ahler metric $g_u$ in a neighborhood of $M$.

\smallskip
\noindent (3) If $u = e^\varphi \hat u $, $\theta = d^c u$ and
$\hat \theta = d^c \hat u$, then one has
$$
\theta = e^\varphi \hat \theta \text{  on   } M.
$$
\endproclaim
\demo{Proof}  Assertion (1) was stated in Theorem 3.4.4 of [CS, p45-46].

The verification of Assertions (2)-(3) is straightforward.
\qed\enddemo

\proclaim{Proposition 3.2} Let $\Omega$ be a bounded, strictly pseudo-convex domain with a smooth
boundary $M = b\Omega$ in a Stein manifold $V^{2n}$, let
$\Cal{O}$ be  the  holomorphic $(1, 0)$-form bundle of $V^{2n}$, and let $K^*$ be the canonical line bundle of $V^{2n}$. Suppose
that  $\dim_{\Bbb R} [ V^{2n} ] = 2n \ge 6$.
Then the following is true.

\smallskip
\noindent
(1) The first Chern class of $\Cal{O}|_\Omega$ is equal to zero, i.e., $c_1(\Cal{O}|_\Omega) = 0$;
Moreover, the first Chern class of canonical line $c_1(K^*|_M) = 0$;

\smallskip
\noindent
(2)  The Ricci curvature form $Ric_g $ of any metric $g = d\theta $ on $\Cal{O}|_\Omega$  is a $d$-exact $(1, 1)$-form on $M$.
Furthermore, $Ric(\xi, \bar \xi)$ is a real number for all $\xi \in T^{(1, 0)}(M)$.
\endproclaim
\demo{Proof} (1) We will use curved version of Kohn-Morrey formula to verify that
$$
c_1(\Cal{O}|_\Omega) = 0. \tag3.2
$$

Recall that the closure $\bar \Omega$ of $\Omega$ is compact. Since $V^{2n}$ is a Stein manifold, there is a
strictly pluri-subharmonic function $\phi_0$. Let $\phi = \lambda \phi_0$ for sufficiently large $\lambda >0 $.
Using Bochner-H\"ormander-Kohn-Morrey formula, we obtain
$$
  H^{(p, q)}(\Omega) = 0, \tag3.3
$$
for all $0 < q < n$, (cf.  Proposition A.4 of [CaWS, p218]).

It is well-known that, for $\dim_{\Bbb C}(\Omega) = n > 2$
$$
H^1(\Omega, \Cal{O}|_\Omega ) = H^{(1, 1)}(\Omega) = 0. \tag3.4
$$

It follows that the first Chern class of $\Cal{O}|_\Omega$ is zero.

Choose a K\"ahler metric $\hat g$ on $\Omega$. Then the Ricci curvature form $\hat \Theta$ is a $d$-exact
$(1, 1)$-form.

The classical Kohn-Rossi theory states that any $\db_b$-closed $(1, 0)$-form on $M = b\Omega$ can be
extend to a unique holomorphic $(1, 0)$-form on the whole $\Omega$. Thus,
$$
          H^{(1, 1)}(M) = 0, \tag3.5
$$
see [KoR].

 It is also known that $c_1(\Cal{O}|_M) = c_1 (K^*|_M) = 0$.

\medskip

 (2)
Let $g_u$ be the K\"ahler metric associated with the K\"ahler form $i \partial \db u $.
The corresponding first Chern curvature form $\Theta_u$ of the K\"ahler metric $g_u$ is a closed
$(1, 1)$-form in a neighborhood of $M$ in $V^{2n}$.

The classical Chern-Weil theory implies that the cohomology class of the
first Chern curvature form $\Theta|_M$  is independent of the choice of the choice of affine connections on $M$.

In fact, if $\theta_u = d^c u$ and $i \partial \db u  > 0$, then $d\theta_u = d d^c u = i \partial \db u  > 0$ gives rise
a K\"ahler metric in a neighborhood of $M$. For any other $\tilde \theta = e^{ 2 \varphi} \theta $, the Ricci curvature
form corresponding to $\tilde \theta$ remains to be of type $(1, 1)$, see Lemma 2.4 of [L2].
\qed\enddemo

We now recall that a result of Lee [L2].
\proclaim{Proposition 3.3} ([L2, Lemma 6.1, p173-174]) Let $M = b\Omega$ and $\Omega \subset V^{2n}$ be as in Main Theorem.
Suppose that $\tilde \theta = e^{2u}\hat \theta$ and $\hat{R}ic$ is the Ricci curvature form corresponding to $\hat \theta$. Then
$\tilde \theta$ is pseudo-Einstein if and only if there is a real solution $u$ satisfying
$$
i \partial_b  \db_b u = \hat{R}ic
$$
\endproclaim
\demo{Proof} By (6.3) of [L2], the trace-less part of
$\tilde{R}ic$ is zero if there is $\varphi$ satisfying
$$
(n+1) i \partial_b  \db_b \varphi = \hat{R}ic
$$
Since $\hat{R}ic$ is a real valued $d$-exact real-valued $(1,
1)$-form by Proposition 3.2 above, we can choose $\varphi$ to be
real-values as well. (Otherwise, let $v = \frac 12 ( \varphi + \bar \varphi)$
instead). \qed\enddemo

\proclaim{Proposition 3.4} Let $M = b\Omega$ and $\Omega \subset
V^{2n}$ be as in Main Theorem. Suppose that the Ricci curvature
form $\hat{R}ic$ form is a $d$-exact $(1, 1)$-form for the contact
1-form $\hat \theta$. Then there always a real-valued function $u$
satisfying
$$
i \partial_b  \db_b u = \hat{R}ic \tag3.6
$$
\endproclaim
\demo{Proof} Choose $\sigma$ such that
$$
d\sigma = \hat{R}ic. \tag3.7
$$
Let $\sigma = \sigma^{(0,1)} + \sigma^{(1, 0)} + \lambda \theta$, where $\sigma^{(0,1)}  $ is $(0, 1)$-component
of $\sigma$. Since $\hat{R}ic$ is of type $(1,1)$, by (3.7) we have
$$
\db_b \sigma^{(0, 1)} = 0. \tag3.8
$$

Because $\dim_{\Bbb C}(\Omega) > 2$, by a Theorem of Kohn that
there is complex-valued function $f$ with
$$
i \db_b f = \sigma^{(0, 1)}, \tag3.9
$$
see [CS, Ch9].

It follows that
$$
i \partial_b  \db_b f = \partial \sigma^{(0, 1)} = (d \sigma)_b = (\hat{R}ic)_b. \tag3.10
$$
Since $(\hat{R}ic)_b$ is real-valued (1,1)-form, choosing $u =
Re\{ f\}$, we are done. \qed\enddemo

We now summarize our result of this section.

\proclaim{Corollary 3.5    } Suppose that $\Omega \subset V^{2n}$ be a compact strictly pseudo-convex domain with smooth boundary in
a Stein manifold $V^{2n}$. Then its boundary $M^{2n-1} = b\Omega$ admits an intrinsic pseudo-Einstein (i.e., partially Einstein)
 metric.
\endproclaim
\demo{Proof} This is a direct consequence of Lemma 3.1 and Propositions 3.2-3.4.
\qed\enddemo

\bigskip
\head 4. Estimates for CR Paneitz operators on $M^3$
\endhead

In the remaining of this paper, we study the so-called CR Paneitz operator
$$
         P_u f = \triangle_b^2 f +  T^2 f + 4 Im \nabla_\beta(A^{\alpha\beta} \nabla_\alpha f), \tag4.1
$$
where $T = J \nabla u $ is the Reeb vector and $A$ is the torsion tensor of the contact form $\theta_u
$.

When the torsion $A$ vanishes,   the formula (4.1) reduces to (1.2).

 It remains
to verify that CR Paneitz operator $P_u$ is a closed operator.

   If $ \hat \theta = e^\varphi \theta_u$ on $M^3$ and $\hat Q $ is the corresponding CR $Q$-curvature of the metric associated with
   the contact form $\hat \theta$, then
   $$
 e^{2 \varphi} \hat Q  = Q + P_u \varphi,
   $$
see  (5.14) of [GG].

Our goal is to show the following result.

\proclaim{Proposition 4.1} Let $\Omega \subset V^{4}$ be an open strictly pseudo-convex
domain with compact closure in a Stein manifold $V^4$ and let $M^{3} = b\Omega$ be its boundary.
Suppose that $g_u$ is the Cheng-Yau Einstein metric on $\Omega$ and $\theta_u(.)  = du(J.)$ is the
corresponding contact 1-form on $M^3$. Then the Paneitz operator $P_u$ is closed:
$$
        \int_{M^3}| P_u f|^2  \ge c  \int_{M^3} |f|^2, \tag4.2
$$
for any real valued function $f \bot ker P_u$, where $c > 0$ is a
constant independent of $f$.
\endproclaim

\bigskip
\noindent
{\bf Remark 4.2:} The constant $c$ in Proposition 4.1 depends mostly on the Tanaka-Webster curvature $R$ and
pseudo-hermitian torsion $A_{11}$ of $(M^3, J, \theta_u)$ respectively. In fact, the following holds:
$$
\int_M  2 (Pf)  f \theta_u \wedge d \theta_u  =  \int_M
[ 3(\triangle_b f)^2  - |Hess_b f |^2  - R |\nabla_b f |^2 - 6 Im\{A_{\overline{11}}
f_1 f_1 \}  ] \theta_u \wedge d \theta_u
$$
where $  \nabla_b$ and $Hess^2_b$ denotes the sub-gradient and sub-Hessian with respect to $(J, \theta_u)$ respectively,
see [CC].

For the proof of Proposition 4.1, we need some notations.

In what follows, we let $\theta = \theta_u$ be the given contact form. The vector $T$ is the characteristic vector
in $T(M)$ such that $\theta(T) = 1$, $(d\theta)(T, .) = 0$.

An (1, 0)-form $\theta^1 \in \Lambda_{(1, 0)}(M^3)$ is called
{\it admissible} if
$$
\theta^1 (T) = 0, d \theta = i h_{1, \bar{1}}  \theta^1 \wedge \theta^{\bar{1}}
$$
for some hermitian metric function $ h_{1, \bar{1}} $.

It is known that
$$
  \Delta_b f = - { f_\alpha}^\alpha - { f_{\bar \alpha}}^{\bar \alpha}
$$
and
$$
\square_b f = 2 (\db^*_b  \db_b + \db_b \db_b^*) f = ( \Delta_b + i T) f = - 2 { f_{\bar \alpha}}^{\bar \alpha}.
$$

Inspired by proof of Proposition 3.4 of [L2], we will express the
CR Paneitz operator $P$ as a product of several closed operators.

     We first consider
     $$
\Cal L  f = d_b^c f + (\Delta_b f) \theta, \tag4.3
     $$
where $\theta$ is the contact $1$-form described above.

\proclaim{Lemma 4.2}  Let $M^3 = b\Omega$, $\theta$, $A$ and $\Cal L $ be as above. Suppose that $\Omega \subset V^4$
 is a strictly pseudo-convex domain in a Stein manifold $V^4$ and that $\Omega$ has compact closure. Then
$\Cal L $ is a closed operator.

Moreover, one has
$$
d[\Cal L  f] = 2 ( {{f_{\bar 1}}^{\bar 1}}_1 + i A_{ 11} f^1 ) \theta \wedge \theta^1 +
2 ( {{f_{ 1}}^{ 1}}_{\bar 1} -  i A_{ \bar{1} \bar{1} } f^{\bar 1} ) \theta \wedge \theta^{\bar 1}.
$$
\endproclaim
\demo{Proof} By Theorem 9.4.2 of [CS], both $d_b^c$ and
$\triangle_b$ are  closed operators for
strictly pseudo-convex compact CR-hypersurfaces. Notice that
$d_b^c f \in [\Lambda_{(1,0)}(M^3) \oplus \Lambda_{(0, 1)}(M^3)]$
is always orthogonal to the 1-form $(\Delta_b f ) \theta$. Hence,
$\Cal L$ is a closed operator.

We will use the proof of Proposition 3.4 of [L2].

The $\theta^1 \wedge \theta^{\bar 1}$-component of $d[\Cal L  f]$ is
$$
i[f_{1 \bar 1} + f_{\bar{1} 1} - (f_{ 1}^{ 1} +  f_{\bar 1 }^{\bar 1}  ) h_{1 \bar 1}]
\theta^1 \wedge \theta^{\bar 1} = 0.
$$

On the other hand, the $\theta  \wedge \theta^{ 1}$-component of $d[\Cal L  f]$ is
$$
[{{f_{ 1}}^{ 1}}_{ 1} + {{f_{\bar 1}}^{\bar 1}}_1 - i f_{1,0} + i A_{11} f^1]
\theta  \wedge \theta^{ 1}. \tag4.4
$$
It is known (cf. [L2, Section 2]) that
$$
{ - {f_{ 1}}^{ 1}}_{ 1} + {{f_{\bar 1}}^{\bar 1}}_1 + i f_{1,0} + i A_{11} f^1  = 0.  \tag4.5
$$
It follows from (4.4) and (4.5) that the $\theta  \wedge \theta^{ 1}$-component of $d[\Cal L  f]$ is equal to
$$
2 ( {{f_{\bar 1}}^{\bar 1}}_1 + i A_{ 11} f^1 ) \theta \wedge \theta^1. \tag4.6
$$

For the same reason, the $\theta  \wedge \theta^{\bar 1}$-component of $d[\Cal L  f]$ is equal to
$$
2 ( {{f_{ 1}}^{ 1}}_{\bar 1} -  i A_{ \bar{1} \bar{1} } f^{\bar 1} ) \theta \wedge \theta^{\bar 1}.
 \tag4.7
$$
This completes the proof.
 \qed\enddemo

 \demo{Proof of Proposition 4.1 } We now consider the composition of operators:
$$
\tilde P f = \partial^*_b [ \big( d (\Cal L f)  \big) \lfloor_{T}
]. \tag4.8
$$
It follows from  that
$$
\big( d (L f)  \big) \lfloor_{T} = 2 ( {{f_{\bar 1}}^{\bar 1}}_1 + i A_{ 11} f^1 )  \theta^1 +
2 ( {{f_{ 1}}^{ 1}}_{\bar 1} -  i A_{ \bar{1} \bar{1} } f^{\bar 1} )  \theta^{\bar 1}. \tag4.9
$$

We observe that $\partial^*_b $ acts on $\Lambda_{(1,0)}(M^3)$
trivially. For real valued function $f$, we further consider
$$
Re [ \tilde P \circ f] = Re [\overline{\square}_b \square_b f] + 4 Im
(A_{\bar{11}}f_1)_1, \tag4.10
$$
where $Re\{z\}$ is the real part of complex number of $z$.

Therefore, it follows from (4.8)-(4.10) that, for real valued
function $f$, we have
$$
        Re [{\tilde P  }  f]
       = \triangle^2_b f + T^2 f + 4 Im (A_{\bar{11}}f_1)_1 = Pf.
       \tag4.11
$$

Thus, the CR Paneitz operator $P$ satisfies
$$
  P f = Re [\tilde P  f ],
                 \tag4.12
$$
where
$$
\tilde P  = \partial^*_b [ \big( d (\Cal L f)  \big) \lfloor_{T}].
$$

A composition of closed operators remains to be a closed operator.

 If $M^3 = b\Omega$ is a compact strictly pseudo-convex
hypersurface in a Stein manifold $V^4$, then $\{ \db_b, d, \partial^*_b,
\Cal L \}$ are closed operators, by Kohn's
$\db_b$-theory (cf. [CS, Theorem 9.4.2, p231]). Theorem 9.4.2 of
[CS] was stated for $\Omega \subset \Bbb C^2$, but its proof is
applicable to $\Omega$ in all Stein manifolds $V^4$  including $\Bbb C^2$.  It is clear that
 the operator $Re$ is a closed operator. Therefore,
 $P = Re{\tilde P}$ is a closed operator as well.
\qed\enddemo
\bigskip

\demo{Proof of Main Theorem}   Main Theorem now follows from Corollary 2.4,  Corollary 3.5 and Proposition 4.1. \qed\enddemo
\bigskip

\bigskip
\noindent
{\bf Acknowledgement.} The first named author would like to thank Professors Alice Chang, Matt Gursky  and Paul Yang for many inspiring conversations.
In particular, authors are grateful to Matt Gursky for his suggestion of using the notion of ``Q-flat metrics".
We are very indebted to Professor Jih-Hsin Cheng for pointing out an overlook on the difference between essentially positive operators
and closed operators in an earlier version of our paper.
Finally, authors would like to thank
the referee for his (or her) many suggestions on re-organizing and expositions of results in this paper,
including the reference section.

\enddocument